\documentclass{article}
\usepackage{graphicx, multicol, latexsym, amsmath, amssymb}
\usepackage[utf8]{inputenc}
\usepackage[export]{adjustbox}
\usepackage{hyperref}
\usepackage{color}
\usepackage{mathtools}
\usepackage[utf8]{inputenc}
\usepackage[english]{babel}
\usepackage{authblk}

\newtheorem{theorem}{Theorem}[section]
\newtheorem{definition}{Definition}[section]
\newtheorem{lemma}[theorem]{Lemma}

\usepackage{geometry} 
\usepackage{booktabs} 
\usepackage{array} 
\usepackage{paralist}
\usepackage{verbatim}

\usepackage{fancyhdr}
\usepackage[nottoc,notlof,notlot]{tocbibind}

\title{A systematic and complete proof of the existence and uniqueness of self-descriptive numbers}
\date{} 
\author{Orazio Sorgonà}
\affil{ \textit{ University of Rome Niccolò Cusano,}\\ \textit{Via Don Carlo Gnocchi, 3, 00166 Rome, Italy}}
\begin{document}  
\maketitle
\begin{abstract}
 A self-descriptive number in a base $b\geq2$ is an integer $n$ of $b$ digits in which the digit $j_i$, $0\leq j_i\leq b-1$ at the position $p_i, \, 0\leq i \leq b-1$, counts how many times the digit $i$ occurs in the number. It’s known that self-descriptive numbers don’t exist the bases $b=2, 3, 6$ and that they  exist and which they are in the bases $b=4, 5$. Also, that at least one defined self-descriptive number exists in each base $b\geq 7$, while it’s unknown if others do, apart from direct negative checks for smaller bases $\geq 7$. All these results, together with a demonstration of the uniqueness for $b \geq 7$, are here obtained through a systematic scheme of proof. The proof is also complete for all the possible cases had been taken into account. \end{abstract}
 \section{Introduction}

 The subject of \textit{self-descriptive numbers} as it will be  defined and discussed in the following has received some interest by mathematicians \cite{libro1}, \cite{libro2},\cite{articolo1},\cite{articolo2}.
 
  A list of these numbers expressed in the base $b=10$ is found on the OEIS \cite{link lista}.\\ 
 
 \indent Part of the mathematical interest in the subject is that although the existence and numerousness of these numbers is very well known for bases $b\leq 6$, and the existence of at least one defined self-descriptive number had been recognized in each base $b\geq 7$, the eventual uniqueness in these greater bases had remained unproven.\\
 
\indent The acknowledged results  seem however to have been obtain'd mostly through direct checks and trial-and-error procedures, easy for smaller bases. On the internet they can be found several amateur algorithms to check the eventual existence of multiple self-descriptive numbers in greater bases.\\
 
\indent In this work instead, the already known results, together with the demonstration of the still unproven uniqueness, are all obtain'd through a systematic scheme of proof that applies the idea of restricted partition of an integer.
\newpage

 \subsection{Self-descriptive numbers}
 \begin{definition}
 A number $n$ of $b$ digits in some base $b\geq 2$, $n=\displaystyle{\sum_{i=0}^{b-1}}\,j_i\,b^{(b-1-i)},\, 0\leq j_i\leq b-1$; represented by the ordered list of its digits $j_i$ at the $b$ positions $p_i,\, i=0,1,\dots,b-1$;\\ is self-descriptive \emph{iff} the digit
 $j_i$ counts how many times the digit $i$ occurs in $n$.
 \end{definition}

Examples
\begin{itemize}
 \item[*] In the base $b=4$ the number $2020$ is self-descriptive, for there are in it :
\begin{itemize}
 \item[-] $2$ instances of the digit “$0$”;
 \item[-] $0$ instances of the digit “$1$”;
 \item[-] and so on
  \end{itemize}

 \item[*] In the base $b=10$ the number $6210001000$ is self-descriptive and unique.
 \end{itemize}
   The following trivially holds
  \begin{lemma}
  In any base $b$ the sum of all the digits of an eventual self-descriptive number is $b$.
  \end{lemma}
  Related to self-descriptive numbers   also \textit{autobiographical numbers} in a base $b$ are considered in literature (\cite{articolo1}). Those are just all the self-descriptive numbers in bases up to $b$ included, expressed in the base $b$.
\subsection{State of the art}
It’s known that \begin{itemize}
 
\item there are no self-descriptive numbers in the  bases $b=2,3,6$; \item in $b=4$ the numbers $2020$ and $1210$ are self-descriptive and unique;
 
\item in $b=5$ the number $21200$ is self-descriptive and unique;
\item in any $b\geq 7$ it’s self-descriptive the number which has similar entries as $6210001000$ , it is the number
\begin{equation}
\label{general}
(b-4)b^{b-1}+(2)b^{b-2}+(1)b^{b-3}+(1)b^3;
\end{equation}
\end{itemize}
 Leaving aside direct checks, proving negative, for smaller  bases $\geq 7$, it has been hitherto unknown if the numbers of the kind of (\ref{general}) are the unique self-descriptive numbers in \emph{any} $b>7$.

  \subsection{Results here proven}
  All the former results are here originally proven in a  systematic way, it is not  through  trial-and-error procedures. Through the same approach it’s also here proven that\\
  
   \emph{A number of the kind of (\ref{general}) is the unique self-descriptive number in each base $b\geq 7$}.

 \section{Theorems and proofs}
 
Given some first entry $J \equiv j_0$, it is the number of instances of “$0$” in $n$,
\[J= (b-m),\quad 1\leq m < b \leftarrow(J=0\quad \mathrm{is}\:\mathrm{ inconsistent} \rightarrow \:m \neq b)\]
  
 there are then  $(b-J-1)=(m-1)$ empty positions left  in the list of the digits, it is $(m-1)$ parts in which some restricted partition of the integer $m$ is to be had.\\
 
  It’s then immediate  that $m\neq1$ for there would be no empty position left.\\

\indent This already implies that
\begin{theorem} \label{teorema due}
There are no self-descriptive numbers in the base $b=2$. 
\end{theorem}

There is then only one kind of restricted partition of $m\geq 2$ in $m-1$ parts, it is 
\begin{equation}\label{part} m=2\: \left(+1 \right)_{\: m-2\:\mathrm{times}\:\mathbf{\mathrm{iff}}\:  m\geq 3.}\end{equation}
  Thus apart from “$J$” at the position $p_0$, and the $J$ instances of “$0$”,\\ all the other digits in a self-descriptive number must be:  \begin{itemize}
  \item[]  “$2$” in $1$ only instance; and
  \item[] “$1$” in eventual instances, yet anyway not more than $2$.
  \end{itemize}

Thus  $ j_1\: = \; 0,1,2 $, and these three and only cases will be now discussed.

 \subsection{Case 1: $j_1=0$}
 There are no instances of ``1'' in the number, thus because of formula (\ref{part}) \begin{itemize}
 \item[] $m=2$, and there is only “$2$” as a nonzero entry at the positions $p_{(i>0)}$.
 \end{itemize}
Consequentially $j_J=2$, as “$J$” must have at least $1$ instance. This means  that “$J$” occurs twice, and, as the only nonzero entries are  “$J$” and “$2$”, so $J=2$. Thus $b=4$, and

  \begin{theorem}\label{teorema quattro primo}
  The number $2020$ is self-descriptive in the base $b=4$.\end{theorem}

\subsection{Case 2: $j_1=1$} Again $(j_J=2\, \rightarrow\, J=2)$ as in case 1.\\

 This  implies

\begin{theorem}\label{teorema cinque}
  The number $21200$ is self-descriptive in the base $b=5$.\end{theorem}
\subsection{Case 3: $j_1=2$}. Then  $j_2=1$, for $p_2$ has a nonzero entry, and it can be only the digit “$1$”; thus

 \paragraph{Subcase 3.1: $J=1$},\\ \indent then $p_J\equiv p_1$ and,  as all the  nonzero entries have already been  considered,  $b=4$ and
 
 \begin{theorem}\label{teorema quattro secondo}
The number $1210$ is self-descriptive in the base $b=4$.
\end{theorem}
 
 \paragraph{Subcase 3.2: $J\neq 1$}, \\ \indent then  $p_J\not\equiv p_1$ and $J \neq 2$ for $j_2=1$ and “$2$” already occurs at the position $p_1$;\\ thus $J=(b-4)$, for $(j_1=2,\;j_2=1,\;j_{J,\, J\neq (1,2)}=1 \rightarrow j_1+j_2+j_J=m=4)$,\\   and $b > 4$, yet $(b\neq5 \leftarrow J\neq 1),\, (b\neq 6 \leftarrow J\neq 2)$.\\ 
 
  Together with the previously deduced  results, this implies
 
\begin{theorem}\label{teorema tre e sei}
There are no self-descriptive numbers in the bases $b=3,\,6$. 
\end{theorem}

 And 
 \begin{theorem}\label{teorema del generale}
  The number $(b-4)b^{b-1}+(2)b^{b-2}+(1)b^{b-3}+(1)b^3$ is self-descriptive in each base $b \geq 7$. \end{theorem}

\subsection{Concluding uniqueness}
  Finally, as $j_1$ cannot assume any other value than $(0,1,2)$ and all the cases had already been considered, it's so proven that
 \begin{theorem}\label{teorema dell'unicità del generale}
The number  $$(b-4)b^{b-1}+(2)b^{b-2}+(1)b^{b-3}+(1)b^3$$
is  unique as self-descriptive number in each base $b\geq7$.
 \end{theorem}
\vfill
\section{Conclusions}
 A systematic scheme of proof that makes a simple use of the idea of restricted partition of an integer has been applied in a complete  proof of the existences and uniquenesses of self-descriptive numbers.\\

 \indent The results are here listed and their proves in the text referred to. 
\begin{itemize}
 \item[(\ref{teorema due},\ref{teorema tre e sei})] There are no self-descriptive numbers in the bases  $b=2,3,6$. 
 \item[(\ref{teorema quattro primo}, \ref{teorema quattro secondo})] In the base $b=4$ there are the two, and only, self-descriptive numbers $2020$ and $1210$.
\item[(\ref{teorema cinque})] In the base $b=5$ there is the one and only self-descriptive number $21200$.
\item[ (\ref{teorema del generale},\ref{teorema dell'unicità del generale})] In each base $b\geq 7$ is self-descriptive and unique  the number 
\[ (b-4)b^{b-1}+(2)b^{b-2}+(1)b^{b-3}+(1)b^3\quad .\]
 \end{itemize}


\newpage
 
 
 
 \begin{thebibliography}{1}
 
\bibitem{libro1} M. Gardner, Mathematical Circus, pp. 128; 135 Prob. 7 Alfred A. Knopf NY 1979.
\bibitem{libro2} Clifford Pickover, Keys to Infinity, Chapter 28, "Chaos in Ontario." New York: Wiley, pp. 217-219, 1995.
\bibitem{articolo1}Tanya Khovanova, A Story of Storytelling Numbers, Math. Horizons, Sep 2009, 14-17.
	\url{https://arxiv.org/abs/0803.0270} [math.CO]
 \bibitem{articolo2}E. Angelini, "Jeux de suites", in Dossier Pour La Science, pp. 32-35, Volume 59 (Jeux math'), April/June 2008, Paris.

\bibitem{link lista}N.J.A. Sloane, Online Encyclopedia of Integer Sequences (OEIS) \url{https://oeis.org/A108551}



 \end{thebibliography}
\end{document}